# ROBUST GLOBAL STABILIZABILITY BY MEANS OF SAMPLED-DATA CONTROL WITH POSITIVE SAMPLING RATE

**Iasson Karafyllis[*] and Costas Kravaris[**]**

[*]**Department of Environmental Engineering, Technical University of Crete, 73100, Chania, Greece**
**email: ikarafyl@enveng.tuc.gr**

[**]**Department of Chemical Engineering, University of Patras**
**1 Karatheodory Str., 26500 Patras, Greece**
**email: kravaris@chemeng.upatras.gr**

**Abstract**
This work proposes a notion of robust reachability of one set from another set under constant control. This notion is used to construct a control strategy, involving sequential set-to-set reachability, which guarantees robust global stabilization of nonlinear sampled data systems with positive sampling rate. Sufficient conditions for robust reachability of one set from another under constant control are also provided. Finally, the proposed method is illustrated through two examples.

# 1. Introduction

Given the finite-dimensional continuous-time system:

$$\dot{x}(t) = f(d(t), x(t), u(t))$$
$$d(t) \in D \subseteq \Re^l,\ x(t) \in \Re^n,\ u(t) \in U \subseteq \Re^m \quad (1.1)$$

where the vector field $f : D \times \Re^n \times U \to \Re^n$ is continuous, $u(t)$ represents the control input and $d(t)$ unknown disturbances or model uncertainty. Consider now a state feedback law $u = k(x)$ to be applied to system (1.1) in discrete time, under zero-order hold with sampling period $h$

$$u(t) = k(x(\tau_i))\ \text{on the interval}\ [\tau_i, \tau_i + h)\ ,\ i = 0,1,2,... \quad (1.2)$$

The resulting closed-loop system is the following hybrid system:

$$\dot{x}(t) = f(d(t), x(t), k(x(\tau_i))),\ t \in [\tau_i, \tau_{i+1})$$
$$\tau_{i+1} = \tau_i + h,\ i = 0,1,2,... \quad (1.3)$$

and the question is how to select the state feedback function $k(x)$ for desirable stability characteristics of (1.3).

There is a large body of literature concerning the above very important and very challenging problem of designing sampled-data feedback stabilizers. In particular, the following lines of attack have been pursued to derive stability results (see also the detailed discussion in review paper [25]):

∗ making use of numerical approximations of the solution of the open-loop system (e.g., in the work of D. Nesic, A. Teel and others, see [6,7,15,19,20,26-32,36]). The results obtained in this way lead to a systematic procedure for the construction of practical, semi-global feedback stabilizers and provide a list of possible reasons that explain the

occasional failure of sampled-data control mechanisms. Recent research takes into account performance and robustness issues as well (see [6,15,17,31]).
* exploiting special characteristics of the system such as homogeneity (see [5]), global Lipschitz conditions (see [8]) or linear structure with uncertainties (see [2] as well as the textbook [33]).
* making use of Linear Matrix Inequalities in the context of hybrid systems (see [9,10,22,35]).
* considering the closed-loop system as a discrete-time system (see for instance [1,23] and the paper [24] which establishes a unified representation for sampled-data systems and discrete-time systems with analytic dynamics). Recent work has established results that characterize the inter-sample behavior of the solutions based on the behavior of the solution of the discrete-time system (see [26]).
* considering the closed-loop system under zero order-hold as a time-delay system. This approach was recently explored in the context of linear systems theory (see [3,4]).

It is important to note that the above very important research results do not provide conditions for <u>global</u> Asymptotic Stability or Input-to-Output Stability for <u>general</u> nonlinear sampled-data systems (usually only semi-global practical stability properties are established or global stability for limited classes of systems).

The goal of the present this work is the development of a design methodology that guarantees robust global asymptotic stability for system (1.3). In this direction, our proposed line of thinking has been motivated by the following considerations:

1) Because the system description (1.1) is in continuous time, it is natural to try to perform controller design in continuous time as a first step, and subsequently, emulate the continuous-time control law under sample-and-hold discretization. But then, the properties of the closed-loop system under continuous-time implementation will not necessarily hold for the closed-loop system under sampled-data control. For example, global closed-loop asymptotic stability in continuous time will <u>not</u>, in general, be preserved under the emulation controller: in general, stability will become local and the size of the domain of attraction will depend on the sampling period (the smaller the sampling period, the larger the domain of attraction). Therefore, a key design issue involves quantifying the capabilities of the emulation controller in terms of size of domain of attraction versus size of sampling period, in order to be able to test whether a satisfactory solution can be obtained in a particular application.

2) If the foregoing emulation design <u>or any other</u> sampled-data controller design guarantees a domain of attraction that is too small in size for a particular application, the immediate question that arises concerns the possibility of extending the control strategy for the purpose of enlargement of the domain of attraction. When the system's initial condition is outside the guaranteed domain of attraction under a given controller, is it possible to find a strategy that can bring the system inside? Assuming that the guaranteed domain of attraction is a positively invariant set, the question is whether it can become reachable from any initial condition outside it – once the system's state is somehow brought inside, the controller will take care of it afterwards.

To be able to address the latter question, some intuitive considerations would be helpful, before a mathematical formulation is developed. In order to drive the system's state to the given target set, the simplest choice of control input that could be tried is constant control; in most applications, such a choice is meaningful on physical grounds as well, and with good chances of success when the target set is large enough. If constant control cannot take the system inside the target set, it will still be able to take it somewhere else, where, on physical grounds, the system will perform satisfactorily. From there, another constant value of the control input can be tried out and, if it still does not hit the target, still another constant control input, … , up until the target set is reached.

The present work will provide a mathematical formulation of the foregoing intuitive idea of sequential reachability from one region of state space to another, ultimately reaching the target attractor. The goal will be to develop and prove conditions under which this intuitive idea will lead to robust global asymptotic stability for the closed-loop system (Theorem 3.1 and Theorem 3.2). In this direction, a new notion of reachability of one set from another under constant control will be proposed (Definition 2.4) and subsequently, this notion will be utilized to establish the main stability results that involve a chain of reachable sets. Simple sufficient conditions to test reachability of one set from another will also be derived (Lemma 2.7 and Lemma 2.9). Finally, the proposed method will be applied to two illustrative examples (Example 4.1 and Example 4.2).

Whenever applicable, the proposed control method has very desirable features, including that
* it guarantees *global* asymptotic stability for the closed-loop system,
* it guarantees *robustness* to perturbations of the sampling schedule,
* it provides means to *determine the maximum allowable sampling period*,
* is *not limited* to special cases where the solution map is available,
* is *not limited* to special cases where the nonlinear term is homogeneous or globally Lipschitz

No other existing method can guarantee all of the above at the same time.

**Notations** Throughout this paper we adopt the following notations:

* For a vector $x \in \Re^n$ we denote by $|x|$ its usual Euclidean norm and by $x'$ its transpose.

* We say that a non-decreasing continuous function $\gamma : \Re^+ \to \Re^+$ is of class $N$ if $\gamma(0) = 0$. We say that a function $\rho : \Re^+ \to \Re^+$ is positive definite if $\rho(0) = 0$ and $\rho(s) > 0$ for all $s > 0$. For the definitions of the classes $K$ and $K_\infty$, see [16]. By $KL$ we denote the set of all continuous functions $\sigma = \sigma(s,t) : \Re^+ \times \Re^+ \to \Re^+$ with the properties: (i) for each $t \geq 0$ the mapping $\sigma(\cdot, t)$ is of class $K$ ; (ii) for each $s \geq 0$, the mapping $\sigma(s, \cdot)$ is non-increasing with $\lim_{t \to +\infty} \sigma(s,t) = 0$.

* Let $D \subseteq \Re^l$ be a non-empty set. By $L_{loc}^\infty(\Re^+; D)$ we denote the class of all Lebesgue measurable and locally essentially bounded mappings $d : \Re^+ \to D$.

* Let $A \subseteq \Re^n$ be a non-empty set. For every $\varepsilon > 0$ we define the $\varepsilon$-neighborhood of $A$ by $N(A, \varepsilon) := \{ y \in \Re^n : dist(y, A) < \varepsilon \}$, where $dist(y, A) = \inf \{ |y - x| : x \in A \}$.

* For every scalar continuously differentiable function $V : \Re^n \to \Re$, $\nabla V(x)$ denotes the gradient of $V$ at $x \in \Re^n$, i.e., $\nabla V(x) = \left( \frac{\partial V}{\partial x_1}(x), \ldots, \frac{\partial V}{\partial x_n}(x) \right)$. We say that a function $V : \Re^n \to \Re$ is positive definite if $V(x) > 0$ for all $x \neq 0$ and $V(0) = 0$. We say that a continuous function $V : \Re^n \to \Re$ is radially unbounded if the following property holds: "$V(x)$ is bounded if and only if $|x|$ is bounded".

## 2. Main Assumptions and Notions for Sampled-Data Systems

In the present work we study control systems of the form (1.1) under the following hypotheses:

**(H1)** $f(d, x, u)$ is continuous with respect to $(d, x, u) \in D \times \Re^n \times U$ and such that for every bounded $S \subset \Re^n \times U$ there exists constant $L \geq 0$ such that

$$(x - y)'(f(d, x, u) - f(d, y, u)) \leq L |x - y|^2 \quad (2.1)$$
$$\forall (x, u, d) \in S \times D , \forall (y, u, d) \in S \times D$$

Hypothesis (H1) is a standard continuity hypothesis and condition (2.1) is characterized as a "one-sided Lipschitz condition on compact sets" in the literature (see [34]). Notice that we do not assume Lipschitz continuity of the vector field $f(d, x, u)$ with respect to $x \in \Re^n$. It is clear that hypothesis (H1) guarantees that for every $(x_0, d, u) \in \Re^n \times L_{loc}^\infty(\Re^+; D) \times L_{loc}^\infty(\Re^+; U)$, there exists a unique solution $x(t)$ of (1.1) with initial condition $x(0) = x_0$ corresponding to inputs $(d, u) \in L_{loc}^\infty(\Re^+; D) \times L_{loc}^\infty(\Re^+; U)$.

**(H2)** There exist a function $a \in K_\infty$ such that

$$|f(d, x, u)| \leq a(|x| + |u|), \quad \forall (u, d, x) \in U \times D \times \Re^n \quad (2.2)$$

Hypothesis (H2) guarantees that $0 \in \Re^n$ is an equilibrium point for (1.1) and is automatically satisfied if $D \subset \Re^l$ is compact.

We next provide a definition of robust global stabilizability of (1.1) by means of bounded sampled-data control with positive sampling rate, which will be used in subsequent developments.

**Definition 2.1:** *We say that the equilibrium point $0 \in \Re^n$ of (1.1) under hypotheses (H1-2) is **robustly globally stabilizable by means of sampled-data control with positive sampling rate**, if there exists a locally bounded mapping $k : \Re^n \to U \subseteq \Re^m$ with $k(0) = 0$ (the feedback function), a function $\kappa \in K_\infty$ with $|f(d, z, k(x))| \leq \kappa(|z| + |x|)$ for all* $(d, z, x) \in D \times \Re^n \times \Re^n$, a constant $h > 0$ (the maximum allowable sampling period) and a function $\sigma \in KL$ such that the following estimate holds for all $(x_0, d, \tilde{d}) \in \Re^n \times L_{loc}^\infty(\Re^+; D) \times L_{loc}^\infty(\Re^+; \Re^+)$ and $t \geq 0$:

$$|x(t)| \leq \sigma(|x_0|, t) \qquad (2.3)$$

where $x(t)$ denotes the solution of the system:

$$\begin{aligned} \dot{x}(t) &= f(d(t), x(t), k(x(\tau_i))), \quad t \in [\tau_i, \tau_{i+1}) \\ \tau_0 &= 0, \tau_{i+1} = \tau_i + h \exp(-\tilde{d}(\tau_i)), i = 0, 1, \ldots \end{aligned} \qquad (2.4)$$

with initial condition $x(0) = x_0$.

We say that the equilibrium point $0 \in \Re^n$ of (1.1) under hypotheses (H1-2) is **robustly globally stabilizable by means of bounded sampled-data control with positive sampling rate**, if the feedback function $k: \Re^n \to U \subseteq \Re^m$ is bounded.

**Remark 2.2:** (a) In this work, the closed-loop system (2.4) will be regarded as a hybrid system, that produces for each $x_0 \in \Re^n$ and for each pair of measurable and locally bounded inputs $d: \Re^+ \to D$, $\tilde{d}: \Re^+ \to \Re^+$, the absolutely continuous function $t \to x(t) \in \Re^n$, produced by the following algorithm:

Step $i$:
1) Given $\tau_i$, calculate $\tau_{i+1}$ using the equation $\tau_{i+1} = \tau_i + h \exp(-\tilde{d}(\tau_i))$,
2) Compute the state trajectory $x(t)$, $t \in [\tau_i, \tau_{i+1})$ as the solution of the differential equation $\dot{x}(t) = f(d(t), x(t), k(x(\tau_i)))$,
3) Calculate $x(\tau_{i+1})$ using the equation $x(\tau_{i+1}) = \lim_{t \to \tau_{i+1}^-} x(t)$.

Hybrid systems of the form (2.4) were considered in [12,13]. Particularly, it was shown that under hypotheses (H1-2) and the hypotheses of Definition 2.1, the hybrid system (2.4) is an autonomous system which satisfies weak semigroup property, the "Boundedness Implies Continuation" property and for which $0 \in \Re^n$ is a robust equilibrium point for system (2.4) from the input $\tilde{d} \in L_{loc}^\infty(\Re^+; \Re^+)$ (see [12]). Moreover, the existence of a function $\sigma \in KL$ that satisfies (2.3) is equivalent to requiring Uniform Robust Global Asymptotic Stability for the closed-loop system (2.4).

(b) Under hypothesis (H2) and the assumption that $k: \Re^n \to U \subseteq \Re^m$ is a locally bounded mapping with $k(0) = 0$, the assumption that there exists $\kappa \in K_\infty$ with $|f(d, z, k(x))| \leq \kappa(|z| + |x|)$ for all $(d, z, x) \in D \times \Re^n \times \Re^n$ is automatically satisfied if the mapping $k: \Re^n \to U \subseteq \Re^m$ is continuous at $x = 0$.

**Remark 2.3:** The reader should notice that the sampling period is allowed to be time-varying. The factor $\exp(-\tilde{d}(\tau_i)) \leq 1$, with $\tilde{d}(t) \geq 0$ some non-negative function of time, is an uncertainty factor in the end-point of the sampling interval. Proving robust global stabilizability of (1.1) by bounded sampled-data feedback with positive sampling rate will guarantee stability of the closed-loop system (2.4) for all sampling schedules with $\tau_{i+1} - \tau_i \leq h$ (robustness to perturbations of the sampling schedule). To understand the importance of robustness to perturbations of the sampling schedule, consider the following situation. Suppose that hardware limitations restrict the sampling period to be $1s$. If we manage to design a sampled-data feedback law with $h(x) \equiv r \geq 2s$, then the application of the feedback control will guarantee stability properties for the closed-loop system even if we "miss measurements" or if we have "delayed measurements" (for example, due to improper operation of the sensor). In such a case robustness to perturbations of the sampling schedule becomes critical. The introduction of the factor $\exp(-\tilde{d}(\tau_i)) \leq 1$ is a mathematical way of introducing perturbations to the sampling schedule; however, it is not unique. Other ways of introducing perturbations of the sampling schedule can be considered.

We next propose a notion of reachability of one set from another set for control systems of the form (1.1), which is going to be utilized for the construction of sampled-data feedback stabilizers in the following section.

**Definition 2.4:** *Consider system (1.1) under hypotheses (H1-2) and let $r > 0$ be a constant. A set $A \subseteq \Re^n$ is $r$-robustly reachable from a set $\Omega \subseteq \Re^n$ for system (1.1) with constant control if there exist $v \in U$, $c \geq 0$ and functions $a, b \in N$ with the following property:*

**(Q)** *For every $x_0 \in \Omega$, $d \in L_{loc}^\infty(\Re^+; D)$, there exists $T(d, x_0) \in [0, c + b(|x_0|)]$ the solution of (1.1) with $u(t) \equiv v$ and initial condition $x(0) = x_0$ exists for all $t \in [0, T(d, x_0) + r]$ and satisfies:*

i) $|x(t)| \leq a(|x_0|)$, for all $t \in [0, T(d, x_0) + r]$

ii) $x(t) \in A$, for all $t \in [T(d, x_0), T(d, x_0) + r]$

iii) $x(t) \in \Omega$, for all $t \in [0, T(d, x_0)]$

**Remark 2.5:** It should be emphasized that $r$-robust reachability of a set with constant control is a much stronger property than simple reachability (see [33]):

(a) property (Q)-(ii) requires that the solution remains in the reachable set for at least time $r$ for all possible disturbances,
(b) property (Q)-(i) requires that the solution remains uniformly bounded for all possible disturbances and for initial conditions in a specified compact set of the state space,
(c) property (Q) requires that the time needed in order reach the set $A \subseteq \Re^n$ is uniformly bounded for all possible disturbances and for initial conditions in a specified compact set of the state space.

**Example 2.6:** Consider the simplified Moore-Greitzer model of a jet engine with no stall presented in [18], described by the planar system:

$$\dot{x}_1 = \frac{3}{2} x_1^2 - \frac{1}{2} x_1^3 + x_2$$
$$\dot{x}_2 = u \qquad (2.5)$$
$$x = (x_1, x_2)' \in \Re^2, u \in \Re$$

The sampled-data stabilizability properties of the jet engine system were studied in [32], where it was shown that system (2.5) can be practically, semiglobally stabilized by sampled-data control with positive sampling rate. Here we study the perturbed version of the jet engine system, i.e., the system:

$$\dot{x}_1 = d_1(t) x_1 + \frac{3}{2} d_2(t) x_1^2 - \frac{1}{2} x_1^3 + x_2$$
$$\dot{x}_2 = u \qquad (2.6)$$
$$x = (x_1, x_2)' \in \Re^2, u \in \Re, d = (d_1(t), d_2(t)) \in [-1,1]^2$$

In this example, we show that the set $\Omega_2 = \{(x_1, x_2) \in \Re^2 : |x_2| \leq 1\}$ is $r$-robustly reachable from the set $\Omega_3 = \{(x_1, x_2) \in \Re^2 : x_2 \leq -1\}$ and from the set $\Omega_4 = \{(x_1, x_2) \in \Re^2 : x_2 \geq 1\}$ for system (2.6) with constant control and $r = 1$.

To prove reachability of $\Omega_2$ from $\Omega_4$, let $v = -1$ and notice that the solution $x(t)$ of (2.6) with initial condition $x_0 = (x_{10}, x_{20})' \in \Omega_4$ satisfies $x_2(t) = x_{20} - t$ for all $t \geq 0$ such that the solution of (2.6) exists. Moreover, we have:

$$\frac{d}{dt}\left(x_1^2(t)\right) = 2d_1(t) x_1^2(t) + 3 d_2(t) x_1^3(t) - x_1^4(t) + 2 x_1(t) x_2(t) \leq 8 x_1^2(t) + x_2^2(t)$$

The above differential inequality in conjunction with the fact that $x_2(t) = x_{20} - t$ gives $x_1^2(t) \leq \left( x_{10}^2 + \frac{1}{8} \max_{\tau \in [0,t]} (x_{20} - \tau)^2 \right) \exp(8 + 8|x_0|)$ for all $t \in [0, 1 + |x_0|]$. Consequently, the solution of (2.6) exists for all $t \in [0, 1 + |x_0|]$. It follows that:

$$|x(t)| \leq \left( |x_0| + \max_{\tau \in [0,t]} |x_{20} - \tau| \right) \exp(4 + 4|x_0|), \text{ for all } t \in [0, 1 + |x_0|] \tag{2.7}$$

Next we show that property (Q) of Definition 2.4 holds with $c := 0$, $b(s) := s \in N$, $a(s) := 2s \exp(4 + 4s) \in N$ and $T(d, x_0) = x_{20} - 1$. Indeed, we have $x(t) \in \Omega_2$ for all $t \in [T(d, x_0), T(d, x_0) + 1]$, $x(t) \in \Omega_2$ for all $t \in [0, T(d, x_0)]$, where $T(d, x_0) = x_{20} - 1$. Moreover, we have $T(d, x_0) \leq c + b(|x_0|)$, where $c := 0$ and $b(s) := s \in N$. Finally, from (2.7) we also obtain $|x(t)| \leq a(|x_0|)$ for all $t \in [0, T(d, x_0) + 1]$, where $a(s) := 2s \exp(4 + 4s)$.

Similarly, we can prove that the set $\Omega_2 = \{(x_1, x_2) \in \Re^2 : |x_2| \leq 1\}$ is $r$-robustly reachable from the set $\Omega_3 = \{(x_1, x_2) \in \Re^2 : x_2 \leq -1\}$ for system (2.6) with constant control and $r = 1$. Particularly, using the same arguments we can show that property (Q) of Definition 2.4 holds with $v = 1$, $c := 0$, $b(s) := s \in N$, $a(s) := 2s \exp(4 + 4s) \in N$ and $T(d, x_0) = x_{20} + 1$. ◁

The following simple lemma provides sufficient conditions for $r$-robust reachability of sets with constant control. More specifically, given a positively invariant set $\Omega \subseteq \Re^n$ for system (1.1), we present conditions for the construction of an appropriate subset $A \subseteq \Omega$, which is $r$-robustly reachable from $\Omega \subseteq \Re^n$ for system (1.1) with constant control for every $r > 0$. The following lemma will be used in the examples of the present work.

**Lemma 2.7:** *Consider system (1.1) under hypotheses (H1-2) and suppose that there exists a set $\Omega \subseteq \Re^n$, a continuously differentiable function $V : \Re^n \to \Re$ and constants $v \in U$, $R \geq 0$, $\delta > 0$ such that*

$$\sup_{d \in D} \nabla V(x) f(d, x, v) \leq -\delta, \text{ for all } x \in \Omega \text{ with } V(x) \geq R \tag{2.8}$$

*Moreover, suppose that there exist functions $a_1, a_2 \in K_\infty$ and a constant $p > 0$, such that for every $x_0 \in \Omega$, $d \in L_{loc}^\infty(\Re^+; D)$, the solution of (1.1) with $u(t) \equiv v$ and initial condition $x(0) = x_0$ exists for all $t \geq 0$ and satisfies $x(t) \in \Omega$, $a_1(|x(t)|) \leq \exp(pt) a_2(|x_0|)$ for all $t \geq 0$. Then for every $r > 0$, the set $A := \Omega \cap \{x \in \Re^n : V(x) \leq R\}$ is $r$-robustly reachable from $\Omega \subseteq \Re^n$ for system (1.1) with constant control.*

**Proof:** Let $r > 0$. Notice that inequality (2.8) guarantees that the set $A := \Omega \cap \{x \in \Re^n : V(x) \leq R\}$ is positively invariant for system (1.1) with $u(t) \equiv v$. Consequently, if $x_0 \in A$ then $x(t) \in A$ for all $t \geq 0$ and $d \in L_{loc}^\infty(\Re^+; D)$. Let arbitrary $x_0 \in \Omega$, with $V(x_0) > R$, $d \in L_{loc}^\infty(\Re^+; D)$ and consider the solution of (1.1) with $u(t) \equiv v$ and initial condition $x(0) = x_0$. Define the set $\{t \geq 0 : x(t) \notin A\}$. Clearly this set is non-empty (since $0 \in \{t \geq 0 : x(t) \notin A\}$). We next claim that $\sup\{t \geq 0 : x(t) \notin A\} \leq \delta^{-1}(V(x_0) - R)$. Suppose that this is not the case. Then there exists $t > \delta^{-1}(V(x_0) - R)$ with $V(x(t)) > R$. Since $A := \Omega \cap \{x \in \Re^n : V(x) \leq R\}$ is positively invariant for system (1.1) with $u(t) \equiv v$, this implies that $V(x(\tau)) > R$ for all $\tau \in [0, t]$. Consequently, it follows from (2.8) that $\frac{d}{d\tau} V(x(\tau)) \leq -\delta$, a.e. on $[0, t]$. Thus we obtain $V(x(t)) \leq V(x_0) - \delta t$, which combined with the hypothesis $t > \delta^{-1}(V(x_0) - R)$ gives $V(x(t)) \leq R$, a contradiction.

Thus, for every $x_0 \in \Omega$, $d \in L_{loc}^\infty(\Re^+; D)$ there exists time $T(d, x_0) \geq 0$ with $T(d, x_0) \leq \delta^{-1} \max\{0, V(x_0) - R\}$ such that $x(t) \in A$, for all $t \in [T(d, x_0), T(d, x_0) + r]$. Furthermore, inequality $T(d, x_0) \leq \delta^{-1} \max\{0, V(x_0) - R\}$ implies $T(d, x_0) \leq c + b(|x_0|)$, where $b(s) := \delta^{-1}\left(\max_{|x| \leq s} \max\{0, V(x) - R\} - \max\{0, V(0) - R\}\right) \in N$ and $c := \delta^{-1} \max\{0, V(0) - R\} \geq 0$. By virtue of the hypotheses of the lemma, $|x(t)| \leq a(|x_0|)$, for all $t \in [0, T(d, x_0) + r]$, where $a(s) := a_1^{-1}(\exp(p(c+r) + pb(s))a_2(s)) \in N$. Consequently, all requirements of Definition 2.4 hold and the set $A := \Omega \cap \{x \in \Re^n : V(x) \leq R\}$ is $r$-robustly reachable from $\Omega \subseteq \Re^n$ for system (1.1) with constant control. The proof is complete. ◁

The following example illustrates how Lemma 2.7 can be used for the establishment of $r$-robust reachability of sets with constant control.

**Example 2.8:** Consider again the perturbed jet engine system (2.6). In this example, we show that the set $A = \{(x_1, x_2) \in \Re^2 : |x_2| \leq 1, |x_1| \leq 4\}$ is $r$-robustly reachable from the set $\Omega_2 = \{(x_1, x_2) \in \Re^2 : |x_2| \leq 1\}$ for system (2.6) with constant control and $r = 1$. Define $\Omega = \Omega_2$, $V(x) = x_1^2$ and $v = 0$. Notice that the solution $x(t)$ of (2.6) with initial condition $x_0 = (x_{10}, x_{20})' \in \Omega_2$ satisfies $x_2(t) = x_{20} \in [-1,1]$ for all $t \geq 0$ such that the solution of (2.6) exists. Moreover, we have:

$$\frac{d}{dt}\left(x_1^2(t)\right) = 2d_1(t)x_1^2(t) + 3d_2(t)x_1^3(t) - x_1^4(t) + 2x_1(t)x_2(t) \leq 8x_1^2(t) + x_2^2(t)$$

The above differential inequality in conjunction with the fact that $x_2(t) = x_{20}$ gives $x_1^2(t) \leq \left(x_{10}^2 + \frac{1}{8}x_{20}^2\right)\exp(8t)$ for all $t \geq 0$. Consequently, the solution of (2.6) exists for all $t \geq 0$ and satisfies $x(t) \in \Omega$ and

$$|x(t)| \leq 2|x_0|\exp(4t), \text{ for all } t \geq 0 \tag{2.9}$$

Moreover, notice that

$$\sup_{d \in [-1,1]^2} 2d_1 x_1^2 + 3d_2 x_1^3 - x_1^4 + 2x_1 x_2 \leq -7, \text{ for all } x \in \Omega \text{ with } V(x) \geq 16 \tag{2.10}$$

It follows from (2.9), (2.10) that the hypotheses of Lemma 2.7 hold with $a_1(s) := s$, $a_2(s) := 2s$, $p = 4$, $\delta = 7$ and $R := 16$. Consequently the set $A := \Omega \cap \{x \in \Re^n : V(x) \leq 16\}$ is $r$-robustly reachable from $\Omega \subseteq \Re^n$ for system (2.6) with constant control and $r = 1$. Notice that $A = \{(x_1, x_2) \in \Re^2 : |x_2| \leq 1, |x_1| \leq 4\}$. ◁

Finally, we end this section with a result that provides links between $r$-robust reachability of sets with constant control and attractor theory for systems without disturbances. Particularly, we show that for every $\varepsilon, r > 0$ an $\varepsilon$-neighborhood of a compact global attractor is $r$-robustly reachable from $\Re^n$. Consequently, knowledge of the dynamics of a control system under constant input may be used for the construction of $r$-robustly reachable sets.

**Lemma 2.9:** *Let $U \subseteq \Re^m$ with $0 \in U$ and consider the control system:*

$$\begin{aligned}\dot{x}(t) &= f(x(t), u(t)) \\ x(t) &\in \Re^n, u(t) \in U\end{aligned} \tag{2.11}$$

*where $f$ is a locally Lipschitz vector field with $f(0) = 0$. Suppose that there exists $v \in U$ such that the dynamical system (2.11) with $u(t) \equiv v$ has a compact global attractor $A \subset \Re^n$. Then for every $\varepsilon, r > 0$ the $\varepsilon$-neighborhood of $A \subset \Re^n$, $N(A, \varepsilon)$ (see notations) is $r$-robustly reachable from $\Re^n$ for system (2.11) with constant control.*

**Proof:** Let $x(t, x_0)$ denote the solution of (2.11) with $u(t) \equiv v$ and initial condition $x(0) = x_0$. Since $A \subset \Re^n$ is a global attractor, for every $\varepsilon, R > 0$ there exists $T(\varepsilon, R) \geq 0$ such that the following implication holds (see [34]):

$$\text{"if } |x_0| \leq R \text{ then } x(t, x_0) \in N(A, \varepsilon) \text{ for all } t \geq T(\varepsilon, R)\text{"} \tag{2.12}$$

Let $g: \Re^+ \to \Re^+$ be defined by $g(s) := T(\varepsilon, k+1) + (s-k)(T(\varepsilon, k+2) - T(\varepsilon, k+1))$ for $s \in [k, k+1)$ and every non negative integer $k$. Clearly, $g: \Re^+ \to \Re^+$ is continuous with $g(s) \geq \min\{T(\varepsilon, [s]+1), T(\varepsilon, [s]+2)\}$, where $[s]$ is the integer part of $s \geq 0$. Define $c := g(0)$ and $b(s) := \max\{g(y) - g(0) : y \in [0, s]\}$. Clearly, $b \in N$ with $g(s) \leq c + b(s)$ for all $s \geq 0$.

Let $x_0 \in \Re^n$ and consider the solution of (2.11) with $u(t) \equiv v$ and initial condition $x(0) = x_0$. By virtue of implication (2.12), there exists $T(x_0) \geq 0$ such that $x(t, x_0) \in N(A, \varepsilon)$ for all $t \geq T(x_0)$. Moreover, $T(x_0) \leq \min\{T(\varepsilon, [|x_0|]+1), T(\varepsilon, [|x_0|]+2)\}$, where $[|x_0|]$ is the integer part of $|x_0|$ and consequently, we obtain $T(x_0) \leq g(|x_0|) \leq c + b(|x_0|)$. Consequently, requirements (ii), (iii) of Property (Q) of Definition 2.4 hold.

We next show that requirement (i) of Property (Q) of Definition 2.4 holds as well for appropriate $a \in N$. Since $A \subset \Re^n$ is bounded, there exists $M > 0$ such that $N(A, \varepsilon) \subseteq \overline{B}_M$, where $\overline{B}_M$ denotes the closed sphere in $\Re^n$ of radius $M > 0$, centered at $0 \in \Re^n$. Consequently, by virtue of implication (2.12), we obtain for all $s \geq 0$:

$$\sup\{|x(t, x_0)| : t \geq 0, |x_0| \leq s\} \leq \max\left(\sup\{|x(t, x_0)| : t \in [0, T(\varepsilon, s)], |x_0| \leq s\}, \sup\{|x(t, x_0)| : t \geq T(\varepsilon, s), |x_0| \leq s\}\right)$$
$$\leq \max\left(\sup\{|x(t, x_0)| : t \in [0, T(\varepsilon, s)], |x_0| \leq s\}, M\right)$$

By virtue of continuity of the mapping $\Re^+ \times \Re^n \ni (t, x_0) \to |x(t, x_0)| \in \Re^+$ and compactness of the set $\{(t, x_0) \in \Re^+ \times \Re^n : t \in [0, T(\varepsilon, s)], |x_0| \leq s\}$, it follows that $\sup\{|x(t, x_0)| : t \in [0, T(\varepsilon, s)], |x_0| \leq s\} < +\infty$. Therefore, for all $s \geq 0$, it holds that $\sup\{|x(t, x_0)| : t \geq 0, |x_0| \leq s\} < +\infty$. By virtue of Lemma 3.5 in [11] there exist a continuous positive function $\mu : \Re^+ \to (0, +\infty)$ and a function $\tilde{a} \in K_\infty$ such that

$$|x(t, x_0)| \leq \mu(t) \tilde{a}(|x_0|) \text{ for all } (t, x_0) \in \Re^+ \times \Re^n \tag{2.13}$$

Define $a(s) := \tilde{a}(s) \max\{\mu(t) : t \in [0, c + b(s) + r]\}$ and notice that since $T(x_0) \leq c + b(|x_0|)$, it follows from (2.13) that requirement (i) of Property (Q) of Definition 2.4 holds as well. The proof is complete. ◁

## 3. Main Results

Our main results are presented below. Theorem 3.1 is an existence result for bounded sampled-data feedback, while Theorem 3.2 is an existence for locally bounded sampled-data feedback. The reader should notice that both theorems do not guarantee continuity of the sampled-data feedback stabilizer.

**Theorem 3.1:** *Consider system (1.1) under hypotheses (H1-2) and suppose the following:*

**(P1)** *There exist a locally bounded mapping $\tilde{k} : \Re^n \to U \subseteq \Re^m$ with $\tilde{k}(0) = 0$, a bounded open set $\Theta \subseteq \Re^n$ which contains a neighborhood of $0 \in \Re^n$, a function $\gamma \in K_\infty$ with $|f(d, z, \tilde{k}(x))| \leq \gamma(|z| + |x|)$ for all $(d, z, x) \in D \times \Re^n \times \Re^n$, a constant $\tilde{h} > 0$ and a function $\sigma \in KL$ such that the following estimate holds for all $(x_0, d, \overline{d}) \in \Theta \times L_{loc}^\infty(\Re^+; D) \times L_{loc}^\infty(\Re^+; \Re^+)$ and $t \geq 0$:*

$$|x(t)| \leq \sigma(|x_0|, t) \;,\; x(t) \in \Theta \tag{3.1}$$

*where $x(t)$ denotes the solution of the system:*

$$\dot{x}(t) = f(d(t), x(t), \tilde{k}(x(\tau_i))) \;,\; t \in [\tau_i, \tau_{i+1})$$
$$\tau_0 = 0, \tau_{i+1} = \tau_i + \tilde{h} \exp(-\overline{d}(\tau_i)), i = 0, 1, \ldots \tag{3.2}$$

*with initial condition* $x(0) = x_0 \in \Theta$.

**(P2)** *There exist sets* $\Omega_j \subseteq \Re^n$, $j = 1,...,N$ *with* $\Omega_1 = \Theta$, $\underset{j=1,...,N}{\cup} \Omega_j = \Re^n$, *such that for each* $j \in \{2,...,N\}$ *the set* $\underset{i=1}{\overset{j-1}{\cup}} \Omega_i$ *is* $r$-*robustly reachable from* $\Omega_j \subseteq \Re^n$ *for system (1.1) with constant control.*

*Then the equilibrium point* $0 \in \Re^n$ *of (1.1) under hypotheses (H1-2) is robustly globally stabilizable by means of bounded sampled-data control with positive sampling rate.*

**Theorem 3.2:** *Consider system (1.1) under hypotheses (H1-2) and suppose that hypothesis (P1) of Theorem 3.1 holds as well as the following hypothesis:*

**(P3)** *There exists a sequence of sets* $\Omega_j \subseteq \Re^n$, $j = 1,2,...$ *with* $\Omega_1 = \Theta$, $\underset{j=1}{\overset{\infty}{\cup}} \Omega_j = \Re^n$, *such that for each* $j \in \{2,3,...\}$ *the set* $\underset{i=1}{\overset{j-1}{\cup}} \Omega_i$ *is* $r$-*robustly reachable from* $\Omega_j \subseteq \Re^n$ *for system (1.1) with constant control. Moreover, for each compact* $K \subseteq \Re^n$, *there exists* $N \in \{2,3,...\}$ *such that* $K \subseteq \underset{j=1}{\overset{N}{\cup}} \Omega_j$.

*Then the equilibrium point* $0 \in \Re^n$ *of (1.1) under hypotheses (H1-2) is robustly globally stabilizable by means of sampled-data control with positive sampling rate.*

**Remark 3.3:** We next provide a brief discussion of hypotheses (P1), (P2), (P3). Hypothesis (P1) is a local hypothesis, which guarantees the existence of a sampled-data feedback, which "works effectively" in the set $\Theta \subseteq \Re^n$. There are many tools in the literature that can be used for the verification of hypothesis (P1) (see for instance [27,31]). It should be emphasized that sampled-data feedback designed by emulation is expected to satisfy hypothesis (P1) for appropriate set $\Theta \subseteq \Re^n$. On the other hand, hypotheses (P2), (P3) are hypotheses of global nature. All tools presented in previous section can be used in order to show the existence of appropriate sets $\Omega_j \subseteq \Re^n$. It should be emphasized that the role of nonlinearities in the verification of hypotheses (P2), (P3) is essential (contrary to hypothesis (P1), which as a local hypothesis depends heavily on the linearization of system (1.1)).

The following lemma can be used for the verification of hypothesis (P1). Its proof can be found in the Appendix.

**Lemma 3.4:** *Consider system (1.1) under hypotheses (H1-2) and suppose that there exists a locally bounded mapping* $\tilde{k}: \Re^n \to U \subseteq \Re^m$ *with* $\tilde{k}(0) = 0$, *a function* $\gamma \in K_\infty$ *with* $\left|f(d,z,\tilde{k}(x))\right| \leq \gamma(|z|+|x|)$ *for all* $(d,z,x) \in D \times \Re^n \times \Re^n$, *a continuous positive function* $\rho: \Re^+ \to \Re^+$, *a positive definite, continuously differentiable and radially unbounded function* $V: \Re^n \to \Re^+$, *constants* $R, \gamma > 0$, $M, L \geq 0$ *such that the following inequalities hold for all* $z \in \Theta$, $x \in \Theta$ *and* $d \in D$:

$$(z-x)' f(d,z,\tilde{k}(x)) \leq L|z-x|^2 + \gamma |x|^2 \tag{3.3}$$

$$\sup\{\nabla V(z) f(d,z,\tilde{k}(x)): d \in D, M|z-x| \leq |z|\} \leq -\rho(V(z)) \tag{3.4}$$

*where* $\Theta := \{x \in \Re^n : V(x) < R\}$. *Let* $x(t)$ *denote the solution of (3.2), initial condition* $x(0) = x_0 \in \Re^n$ *and corresponding to input* $(d, \tilde{d}) \in L_{loc}^\infty(\Re^+; D) \times L_{loc}^\infty(\Re^+; \Re^+)$. *Then hypothesis (P1) of Theorem 3.1 holds with* sampling period $\tilde{h} > 0$ satisfying $\tilde{h} < \frac{1}{2L} \ln\left(1 + \frac{L}{\gamma}\left(\frac{1}{1+M}\right)^2\right)$, for the case $L > 0$ or $\tilde{h} < \frac{1}{2\gamma}\left(\frac{1}{1+M}\right)^2$, for the case $L = 0$.

The rest of the section is devoted to the proof of Theorems 3.1, 3.2.

**Proof of Theorem 3.1:** Define recursively the following sets by the following formulae:

$$C_i = \Omega_i \setminus B_{i-1}, \; B_i = B_{i-1} \cup \Omega_i, \; i = 2,\ldots,N \tag{3.5a}$$

with

$$C_1 = \Omega_1 = \Theta, \; B_1 = \Omega_1 = \Theta \tag{3.5b}$$

Notice that $B_i = \bigcup_{k=i,\ldots,i} \Omega_k = \bigcup_{k=i,\ldots,i} C_k$ for all $i = 1,\ldots,N$. Consequently, by virtue of hypothesis (P2), $B_N = \Re^n$. Let $v_i \in U$ be the constant control that guarantees property (Q) of Definition 2.4 for every set $\Omega_i$ with $i > 1$. We define:

$$k(x) = v_i \text{ if } x \in C_i \text{ with } i > 1 \tag{3.6a}$$

$$k(x) = \tilde{k}(x) \text{ if } x \in C_1 = \Theta \tag{3.6b}$$

$$h = \min\{\tilde{h}, r\} \tag{3.6c}$$

Notice that since $\tilde{k} : \Re^n \to U \subseteq \Re^m$ is locally bounded, it follows that the mapping $k : \Re^n \to U$ as defined by (3.6a,b) is bounded.

We next claim that there exists a function $\kappa \in K_\infty$ with $|f(d,z,k(x))| \leq \kappa(|z| + |x|)$ for all $(d,z,x) \in D \times \Re^n \times \Re^n$ and a function $\sigma \in KL$ such that estimate (2.3) holds for all $(x_0, d, \tilde{d}) \in \Re^n \times L^\infty_{loc}(\Re^+; D) \times L^\infty_{loc}(\Re^+; \Re^+)$ and $t \geq 0$ for the solution $x(t)$ of (2.4) with initial condition $x(0) = x_0$ and corresponding to inputs $(d, \tilde{d}) \in L^\infty_{loc}(\Re^+; D) \times L^\infty_{loc}(\Re^+; \Re^+)$.

Notice that by virtue of hypotheses (P1), (H2) the function $\Re^+ \ni s \to \tilde{\kappa}(s) := \sup\{|f(d,z,k(x))| : |z| + |x| \leq s, d \in D\}$ is a non-decreasing function which satisfies $\tilde{\kappa}(s) \leq \gamma(s)$ for $s \geq 0$ sufficiently small and $|f(d,z,k(x))| \leq \tilde{\kappa}(|z| + |x|)$ for all $(d,z,x) \in D \times \Re^n \times \Re^n$. It turns out that $\tilde{\kappa}$ can be bounded from above by the $K_\infty$ function $\kappa$ defined by $\kappa(s) := s + \frac{1}{s}\int_s^{2s} \tilde{\kappa}(w)dw$ for $s > 0$ and $\kappa(0) = 0$. Consequently, there exists a function $\kappa \in K_\infty$ with $|f(d,z,k(x))| \leq \kappa(|z| + |x|)$ for all $(d,z,x) \in D \times \Re^n \times \Re^n$.

In order to show the existence of a function $\sigma \in KL$ such that estimate (2.3) holds for all $(x_0, d, \tilde{d}) \in \Re^n \times L^\infty_{loc}(\Re^+; D) \times L^\infty_{loc}(\Re^+; \Re^+)$ and $t \geq 0$ for the solution $x(t)$ of (2.4) with initial condition $x(0) = x_0$ and corresponding to inputs $(d, \tilde{d}) \in L^\infty_{loc}(\Re^+; D) \times L^\infty_{loc}(\Re^+; \Re^+)$, we need to show the following things:

- for every $s > 0$, it holds that

$$\sup\{|x(t)| ; t \geq 0, \; |x_0| \leq s, \; (d, \tilde{d}) \in L^\infty_{loc}(\Re^+; D) \times L^\infty_{loc}(\Re^+; \Re^+)\} < +\infty$$
**(Robust Lagrange Stability)**

- for every $\varepsilon > 0$ there exists a $\delta := \delta(\varepsilon) > 0$ such that:

$$\sup\{|x(t)|; t \geq 0, |x_0| \leq \delta, (d,\tilde{d}) \in L_{loc}^{\infty}(\Re^+;D) \times L_{loc}^{\infty}(\Re^+;\Re^+)\} \leq \varepsilon$$
**(Robust Lyapunov Stability)**

- for every $\varepsilon > 0$ and $s \geq 0$, there exists a $\tau := \tau(\varepsilon, s) \geq 0$, such that:

$$\sup\{|x(t)|; t \geq \tau, |x_0| \leq s, (d,\tilde{d}) \in L_{loc}^{\infty}(\Re^+;D) \times L_{loc}^{\infty}(\Re^+;\Re^+)\} \leq \varepsilon$$
**(Uniform Attractivity)**

The above properties guarantee the existence of a function $\sigma \in KL$ such that estimate (2.3) holds for all $(x_0, d, \tilde{d}) \in \Re^n \times L_{loc}^{\infty}(\Re^+;D) \times L_{loc}^{\infty}(\Re^+;\Re^+)$ and $t \geq 0$ for the solution $x(t)$ of (2.4) with initial condition $x(0) = x_0$ and corresponding to inputs $(d, \tilde{d}) \in L_{loc}^{\infty}(\Re^+;D) \times L_{loc}^{\infty}(\Re^+;\Re^+)$. Indeed, we may define $\tilde{\sigma}(s,t) := \sup\{|x(\xi)|; \xi \geq t, |x_0| \leq s, (d,\tilde{d}) \in L_{loc}^{\infty}(\Re^+;D) \times L_{loc}^{\infty}(\Re^+;\Re^+)\}$ for all $s, t \geq 0$. By defining $\tilde{\sigma}(s,t) := \tilde{\sigma}(s,0)$ for all $s \geq 0$, $t \in [-1,0)$, the desired $\sigma \in KL$ can be defined by

$$\sigma(s,t) = s\exp(-t) + \frac{1}{s}\int_{t-1}^{t}\int_{s}^{2s}\tilde{\sigma}(\xi,w)d\xi\,dw \text{ for all } t \geq 0, s > 0 \text{ and } \sigma(0,t) = 0 \text{ for all } t \geq 0.$$

Since the solution of (2.4) with $x(0) = x_0$ corresponding to inputs $(d,\tilde{d}) \in L_{loc}^{\infty}(\Re^+;D) \times L_{loc}^{\infty}(\Re^+;\Re^+)$ coincides with the solution of (3.2) with same initial condition corresponding to inputs $(d,\bar{d}) \in L_{loc}^{\infty}(\Re^+;D) \times L_{loc}^{\infty}(\Re^+;\Re^+)$ with $\bar{d}(t) = \tilde{d}(t) + \ln\left(\frac{\tilde{h}}{h}\right)$, it follows that Robust Lyapunov Stability is an immediate consequence of hypothesis (P1) (notice that $\Theta \subseteq \Re^n$ contains a neighborhood of $0 \in \Re^n$). Thus we are left with the proofs of Robust Lagrange Stability and Uniform Attractivity. Robust Lagrange Stability and Uniform Attractivity will be shown with the help of the following fact, which is shown in the Appendix. Let $c_i \geq 0$ and the functions $a_i, b_i \in N$ that guarantee property (Q) of Definition 2.4 for every non-empty set $\Omega_i$ with $i > 1$ and let $a(s) := \max_{i=2,\ldots,N} a_i(s)$, $b(s) := \max_{i=2,\ldots,N} b_i(s)$, $c := \max_{i=2,\ldots,N} c_i$.

**FACT:** Let $\tilde{d} \in L_{loc}^{\infty}(\Re^+;\Re^+)$ and $\pi(\tilde{d}) := \{\tau_0, \tau_1, \tau_2, \ldots\}$ (the set of sampling times), where $\tau_0 = 0$ and $\tau_{i+1} = \tau_i + h\exp(-\tilde{d}(\tau_i))$ for $i \geq 0$. If $x(\tau_i) \in C_k$ for certain $k > 1$, then for every $d \in L_{loc}^{\infty}(\Re^+;D)$ there exists $\tau \in \pi(\tilde{d}) \cap [\tau_i, \tau_i + c + b(|x(\tau_i)|) + r]$ and $m \in \{1,\ldots,k-1\}$ such that $x(\tau) \in C_m$. Moreover, $|x(t)| \leq a(|x(\tau_i)|)$ for all $t \in [\tau_i, \tau]$.

Since $\bigcup_{k=1,\ldots,N} C_k = \Re^n$ and since $C_1 = \Theta$, the above fact implies that for every $(x_0, d, \tilde{d}) \in \Re^n \times L_{loc}^{\infty}(\Re^+;D) \times L_{loc}^{\infty}(\Re^+;\Re^+)$ the solution $x(t)$ of (2.4) with initial condition $x(0) = x_0$ and corresponding to inputs $(d,\tilde{d}) \in L_{loc}^{\infty}(\Re^+;D) \times L_{loc}^{\infty}(\Re^+;\Re^+)$, satisfies $x(\tau) \in \Theta$ for certain $\tau \in \pi(\tilde{d}) \cap [0, Nc + Nb(a^{(N)}(|x_0|)) + Nr]$ and $|x(t)| \leq a^{(N)}(|x_0|)$ for all $t \in [0,\tau]$, where $a^{(N)} = \underbrace{a \circ a \circ \ldots \circ a}_{N \text{ times}}$. By virtue of (3.1), it follows that $|x(t)| \leq \sigma\left(a^{(N)}(|x_0|), t-\tau\right)$ for all $t \geq \tau$. The properties of the $KL$ functions in conjunction with the previous estimate of the solution imply the Uniform Attractivity property. Moreover, we have $|x(t)| \leq \sigma\left(a^{(N)}(|x_0|), 0\right)$, for all $t \geq 0$ (Uniform Lagrange Stability). The proof is complete. ◁

**Proof of Theorem 3.2:** Define recursively the following sets by the following formulae:

$$C_i = \Omega_i \setminus B_{i-1}, \quad B_i = B_{i-1} \cup \Omega_i, \quad i = 2, 3, \ldots, \tag{3.7a}$$

with

$$C_1 = \Omega_1 = \Theta, \quad B_1 = \Omega_1 = \Theta \tag{3.7b}$$

Notice that $B_i = \bigcup_{k=i,\ldots,i} \Omega_k = \bigcup_{k=i,\ldots,i} C_k$ for all $i = 1,2,3,\ldots$. Let $v_i \in U$ be the constant control that guarantees property (Q) of Definition 2.4 for every set $\Omega_i$ with $i > 1$. We define $k: \Re^n \to U \subseteq \Re^m$ with $k(0) = 0$ and $h > 0$ by (3.6).

Since for each compact $K \subseteq \Re^n$, there exists $N \in \{2,3,\ldots\}$ such that $K \subseteq \bigcup_{j=1}^{N} \Omega_j = \bigcup_{j=1}^{N} C_j$ and $\tilde{k}: \Re^n \to U \subseteq \Re^m$ is locally bounded, it follows that the mapping $k: \Re^n \to U$ as defined by (3.6a,b) is locally bounded.

The proof of the existence of $\kappa \in K_\infty$ with $|f(d, z, k(x))| \leq \kappa(|z| + |x|)$ for all $(d, z, x) \in D \times \Re^n \times \Re^n$ follows exactly the same procedure with the proof of Theorem 3.1. Moreover, as in the proof of Theorem 3.1, in order to complete the proof we need to show Robust Lagrange Stability, Robust Lyapunov Stability and Uniform Attractivity for system (2.4).

The proof of Robust Lyapunov Stability follows exactly the same arguments as in the proof of Theorem 3.1. Thus we are left with the proofs of Robust Lagrange Stability and Uniform Attractivity.

Let $s \geq 0$ and consider the closed ball $\{x \in \Re^n : |x| \leq s\}$. By virtue of hypothesis (P3) there exists $N \in \{2,3,\ldots\}$ such that $\{x \in \Re^n : |x| \leq s\} \subseteq \bigcup_{j=1}^{N} \Omega_j$. Let $c_i \geq 0$ and the functions $a_i, b_i \in N$ that guarantee property (Q) of Definition 2.4 for every $\Omega_i$ with $i > 1$ and let $a(s) := \max_{i=2,\ldots,N} a_i(s)$, $b(s) := \max_{i=2,\ldots,N} b_i(s)$, $c := \max_{i=2,\ldots,N} c_i$. The following fact is a simple extension of the fact used in the proof of Theorem 3.1.

**FACT:** Let $\tilde{d} \in L^\infty_{loc}(\Re^+; \Re^+)$ and $\pi(\tilde{d}) := \{\tau_0, \tau_1, \tau_2, \ldots\}$ (the set of sampling times), where $\tau_0 = 0$ and $\tau_{i+1} = \tau_i + h \exp(-\tilde{d}(\tau_i))$ for $i \geq 0$. If $x(\tau_i) \in C_k$ for certain $k \in \{2, \ldots, N\}$, then for every $d \in L^\infty_{loc}(\Re^+; D)$ there exists $\tau \in \pi(\tilde{d}) \cap [\tau_i, \tau_i + c + b(|x(\tau_i)|) + r]$ and $m \in \{1, \ldots, k-1\}$ such that $x(\tau) \in C_m$. Moreover, $|x(t)| \leq a(|x(\tau_i)|)$ for all $t \in [\tau_i, \tau]$.

Since $\{x \in \Re^n : |x| \leq s\} \subseteq \bigcup_{k=1,\ldots,N} C_k$ and since $C_1 = \Theta$, the above fact implies that for every $(x_0, d, \tilde{d}) \in \{x \in \Re^n : |x| \leq s\} \times L^\infty_{loc}(\Re^+; D) \times L^\infty_{loc}(\Re^+; \Re^+)$ the solution $x(t)$ of (2.4) with initial condition $x(0) = x_0$ and corresponding to (arbitrary) inputs $(d, \tilde{d}) \in L^\infty_{loc}(\Re^+; D) \times L^\infty_{loc}(\Re^+; \Re^+)$, satisfies $x(\tau) \in \Theta$ for certain $\tau \in \pi(\tilde{d}) \cap [0, Nc + Nb(a^{(N)}(|x_0|)) + Nr]$ and $|x(t)| \leq a^{(N)}(|x_0|)$ for all $t \in [0, \tau]$, where $a^{(N)} = \underbrace{a \circ a \circ \ldots \circ a}_{N \text{ times}}$. By virtue of (3.1), it follows that $|x(t)| \leq \sigma(a^{(N)}(|x_0|), t - \tau)$ for all $t \geq \tau$. The properties of the $KL$ functions in conjunction with the previous estimate of the solution imply the Uniform Attractivity property. Moreover, we have $|x(t)| \leq \sigma(a^{(N)}(|x_0|), 0)$, for all $t \geq 0$ (Uniform Lagrange Stability). The proof is complete. ◁

## 4. Examples and Applications

In this section examples are presented, which illustrate how the main results of the present work (Theorem 3.1 and Theorem 3.2) can be used for the construction of robust sampled-data feedback stabilizers.

**Example 4.1:** We consider again the perturbed jet engine system (2.6). Here, we intend to prove that the perturbed jet engine system (2.6) is robustly globally stabilizable by means of bounded sampled-data control with positive sampling rate. The proof will exploit Theorem 3.1.

Consider the function

$$V(x) = \frac{1}{2}x_1^2 + \frac{1}{2}(x_2 + 5x_1)^2 \qquad (4.1)$$

Notice that the set $A = \{(x_1, x_2) \in \Re^2 : |x_2| \leq 1, |x_1| \leq 4\}$ is a subset of $\Theta = \{x \in \Re^2 : V(x) < \frac{457}{2} + \varepsilon\}$ for all $\varepsilon > 0$. Consequently, Examples 2.6 and 2.8 show that the sets $\Omega_4 = \{(x_1, x_2) \in \Re^2 : x_2 \geq 1\}$, $\Omega_3 = \{(x_1, x_2) \in \Re^2 : x_2 \leq -1\}$, $\Omega_2 = \{(x_1, x_2) \in \Re^2 : |x_2| \leq 1\}$, $\Omega_1 = \Theta$, satisfy hypothesis (P2) of Theorem 3.1. We next show that hypothesis (P1) of Theorem 3.1 is also satisfied for the function $V$ defined by (4.1).

The derivative of $V$ along the trajectories of system (2.6) is expressed by the following equality for all $(d, x) \in [-1,1]^2 \times \Re^2$:

$$\nabla V(x) \begin{bmatrix} d_1 x_1 + \frac{3}{2} d_2 x_1^2 - \frac{1}{2} x_1^3 + x_2 \\ u \end{bmatrix} = x_1 \left( d_1 x_1 + \frac{3}{2} d_2 x_1^2 - \frac{1}{2} x_1^3 + x_2 \right) + (x_2 + 5x_1) \left( u + 5d_1 x_1 + \frac{15}{2} d_2 x_1^2 - \frac{5}{2} x_1^3 + 5x_2 \right)$$

Using the inequalities $\frac{3}{2} d_2 x_1^3 \leq \frac{1}{4} x_1^4 + \frac{9}{4} x_1^2$, $5|x_1||x_2 + 5x_1| \leq \frac{1}{4} x_1^2 + 25(x_2 + 5x_1)^2$, $\frac{15}{2} x_1^2 |x_2 + 5x_1| \leq \frac{1}{4} x_1^4 + \frac{225}{4}(x_2 + 5x_1)^2$, we obtain the following inequality for all $(d, x) \in [-1,1]^2 \times \Re^2$:

$$\nabla V(x) \begin{bmatrix} d_1 x_1 + \frac{3}{2} d_2 x_1^2 - \frac{1}{2} x_1^3 + x_2 \\ u \end{bmatrix} \leq -\frac{3}{2} x_1^2 - \frac{11}{4}(x_2 + 5x_1)^2 + (x_2 + 5x_1) \left( u + 421 x_1 - \frac{5}{2} x_1^3 + 89 x_2 \right) \qquad (4.2)$$

We define $\tilde{k}(x) = -421 x_1 - 89 x_2 + \frac{5}{2} x_1^3$ and we obtain from (4.2) for all $(d, x, z) \in [-1,1]^2 \times \Re^2 \times \Re^2$:

$$\nabla V(z) \begin{bmatrix} d_1 z_1 + \frac{3}{2} d_2 z_1^2 - \frac{1}{2} z_1^3 + z_2 \\ \tilde{k}(x) \end{bmatrix} \leq -\frac{3}{2} z_1^2 - 2(z_2 + 5z_1)^2 + \left( 421^2 + \frac{225}{16}(z_1^2 + x_1^2)^2 \right) |z_1 - x_1|^2 + 89^2 |z_2 - x_2|^2$$

If $z \in \Theta$, $x \in \Theta$, where $\Theta := \{x \in \Re^n : V(x) < R\}$ (which gives $x_1^2 \leq 2R$, $z_1^2 \leq 2R$), the above inequality implies:

$$\nabla V(z) \begin{bmatrix} d_1 z_1 + \frac{3}{2} d_2 z_1^2 - \frac{1}{2} z_1^3 + z_2 \\ \tilde{k}(x) \end{bmatrix} \leq -\frac{3}{2} z_1^2 - 2(z_2 + 5z_1)^2 + (421^2 + 225 R^2) |z - x|^2 \qquad (4.3)$$

Notice that $z_1^2 + \frac{3}{2}(z_2 + 5z_1)^2 \geq \frac{3}{80} |z|^2$. Consequently, if $\sqrt{\frac{80}{3}(421^2 + 225 R^2)} |z - x| \leq |z|$, $z \in \Theta$, $x \in \Theta$, we obtain from (4.3) for all $d \in [-1,1]^2$:

$$\nabla V(z) \begin{bmatrix} d_1 z_1 + \frac{3}{2} d_2 z_1^2 - \frac{1}{2} z_1^3 + z_2 \\ \tilde{k}(x) \end{bmatrix} \leq -\frac{1}{2} z_1^2 - \frac{1}{2}(z_2 + 5z_1)^2 = -V(z) \qquad (4.4)$$

Inequality (4.4) shows that inequality (3.2) holds with $R = \frac{457}{2} + \varepsilon$, $M = \sqrt{\frac{80}{3}(421^2 + 225 R^2)}$ and $\rho(s) := s$, for all $\varepsilon > 0$. Moreover, we obtain for all $z \in \Theta$, $x \in \Theta$, $d \in [-1,1]^2$:

$$(z-x)'f(d,z,\tilde{k}(x)) = (z_1 - x_1)\left(d_1 z_1 + \frac{3}{2}d_2 z_1^2 - \frac{1}{2}z_1^3 + z_2\right) + (z_2 - x_2)\tilde{k}(x)$$
$$\leq \left(\frac{7}{2} + \sqrt{2R}\right)|z-x|^2 + \frac{1}{2}|x|^2 + \left(\frac{9R}{4} + \frac{R^2}{2}\right)x_1^2 + \frac{1}{2}(\tilde{k}(x))^2 \qquad (4.5)$$

Notice that $|\tilde{k}(x)| \leq 89|x_2| + (5R-421)|x_1| \leq (5R-332)|x|$, for all $x \in \Theta$. Consequently, inequality (4.5) shows that inequality (3.1) holds with $L = \frac{7}{2} + \sqrt{2R}$ and $\gamma = \frac{9R}{4} + \frac{R^2+1}{2} + \frac{1}{2}(5R-332)^2$.

It follows from Theorem 3.1 that the perturbed jet engine system (2.6) is robustly globally stabilizable by means of bounded sampled-data control with positive sampling rate. Since Theorem 3.1 is proved constructively, a bounded sampled-data feedback can be suggested. Particularly, following the proof of Theorem 3.1, the following discontinuous feedback law:

$$k(x) = -421 x_1 - 89 x_2 + \frac{5}{2} x_1^3, \text{ if } x \in C_1 = \left\{x \in \Re^2 : V(x) < \frac{457}{2} + \varepsilon\right\},$$

$$k(x) = 0, \text{ if } x \in C_2 = \left\{(x_1, x_2) \in \Re^2 : |x_2| \leq 1, V(x) \geq \frac{457}{2} + \varepsilon\right\},$$

$$k(x) = 1, \text{ if } x \in C_3 = \left\{(x_1, x_2) \in \Re^2 : x_2 < -1, V(x) \geq \frac{457}{2} + \varepsilon\right\},$$

$$k(x) = -1, \text{ if } x \in C_1 = \left\{(x_1, x_2) \in \Re^2 : x_2 > 1, V(x) \geq \frac{457}{2} + \varepsilon\right\}$$

is a robust sampled-data feedback stabilizer for system (2.6) for all $\varepsilon > 0$. In Figures 1, 2, 3 it is shown the evolution of the states for the closed-loop system (2.6) with

$$u(t) = k(x(\tau_i)) \quad, \quad t \in [\tau_i, \tau_{i+1})$$
$$\tau_0 = 0, \tau_{i+1} = \tau_i + h\exp(-\tilde{d}(\tau_i)), i = 0,1,... \qquad (4.6)$$

The parameters $h, \varepsilon$ were selected to be $h = \varepsilon = 0.001$ and the initial state is $x_1(0) = 10$, $x_2(0) = 2$. In Figure 1 it is shown the evolution of the states of the closed-loop system (2.6) with (4.6) corresponding to inputs $d_1(t) = \tilde{d}(t) \equiv 0$, $d_2(t) \equiv 1$ and in Figure 2 it is shown the evolution of the states of the closed-loop system (2.6) with (4.6) corresponding to inputs $d_1(t) \equiv 1$, $d_2(t) = \sin(t)$, $\tilde{d}(t) \equiv 0$. Finally, in Figure 3 we tested the performance of the system to additional perturbations of the sampling schedule: the inputs were selected $d_1(t) \equiv 1$, $d_2(t) \equiv 1$, $\tilde{d}(t) = |\sin(t)|$. It is clear that in all cases the closed-loop system presents fast convergence of the states to the equilibrium point. The sampling time $\tau_i$ ($i = 0,1,2,...$) that the state trajectory enters the set $\Theta = \left\{x \in \Re^2 : V(x) < \frac{457}{2} + \varepsilon\right\}$ is the time where the derivative $\dot{x}_2(t)$ presents an abrupt jump (from the value -1 to a negative value with large absolute value).

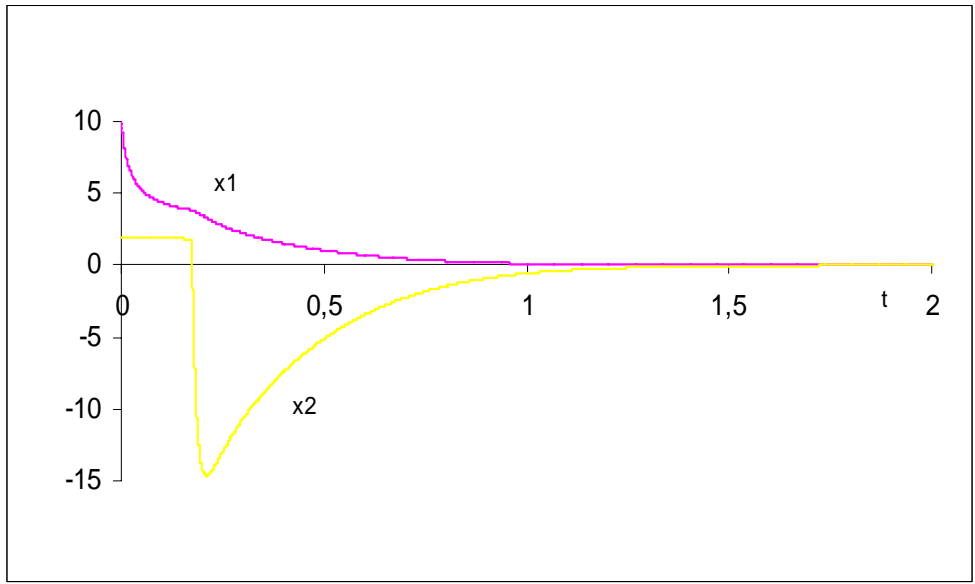

**Figure 1:** The evolution of the states of the closed-loop system (2.6) with (4.6) corresponding to inputs $d_1(t) = \tilde{d}(t) \equiv 0$, $d_2(t) \equiv 1$

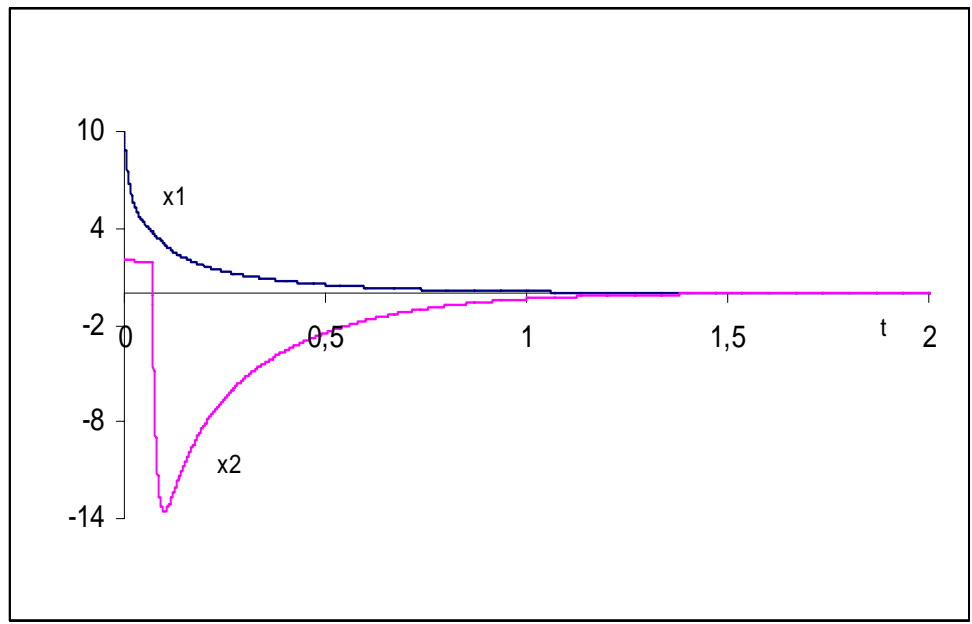

**Figure 2:** The evolution of the states of the closed-loop system (2.6) with (4.6) corresponding to inputs $d_1(t) \equiv 1$, $d_2(t) = \sin(t)$, $\tilde{d}(t) \equiv 0$

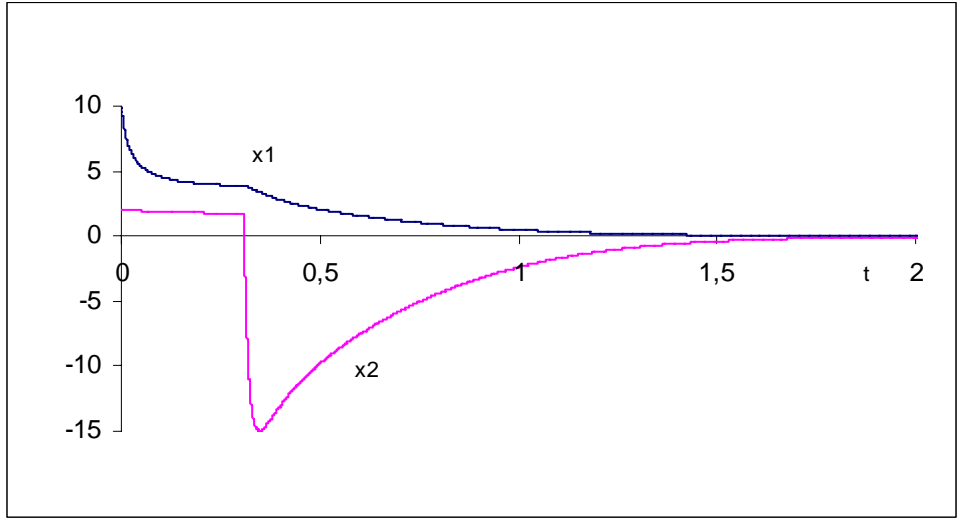

**Figure 3:** The evolution of the states of the closed-loop system (2.6) with (4.6) corresponding to inputs $d_1(t) \equiv 1$, $d_2(t) \equiv 1$, $\tilde{d}(t) = |\sin(t)|$

It should be emphasized that other selections for feedback can be constructed (using different control Lyapunov functions from the quadratic one that we used in this work).  ◁

**Example 4.2:** This example illustrates the use of Theorem 3.2 for the construction of a globally stabilizing sampled-data feedback. Consider the scalar system:

$$\dot{x} = a(x) + u$$
$$x \in \Re, u \in (-\infty, 0] \qquad (4.7)$$

where $a : \Re \to \Re$ is a locally Lipschitz function with $a(0) = 0$ and $a(x) > 0$ for all $x \neq 0$. We claim that system (4.7) satisfies hypotheses (P1), (P3) and consequently, (by virtue of Theorem 3.2) it is robustly globally stabilizable by means of sampled-data control with positive sampling rate.

In order to show the validity of hypothesis (P1), define $\Theta = (-\infty, 2)$, $\tilde{h} = \dfrac{1}{L+1}$ and

$$\tilde{k}(x) := \begin{cases} 0 & \text{for } x \in (-\infty, 0] \\ -(L+1)x & \text{for } x > 0 \end{cases} \qquad (4.8)$$

where $L > 0$ is the Lipschitz constant that satisfies

$$a(x) \leq Lx, \quad \forall x \in [0, 2] \qquad (4.9)$$

The solution of $\dot{x} = a(x) + \tilde{k}(x_0)$ starting at $x(0) = x_0 \in (0, 2)$ satisfies $x(t) \leq x_0$ as long as the solution exists, since by virtue of (4.9) we have $a(x_0) - (L+1)x_0 \leq -x_0 < 0$. Moreover, as long as the solution satisfies $x(t) \geq 0$, it holds that $-(L+1)x_0 \leq \dot{x} \leq -x_0$, which directly implies $(1 - (L+1)t)x_0 \leq x(t) \leq (1-t)x_0$. A simple contradiction argument shows that $0 \leq x(t) \leq (1-t)x_0$ for all $t \in [0, \tilde{h}]$. Consequently, $x(t) \in \Theta$ and $|x(t)| \leq \exp(-t)|x_0|$, for all $t \in [0, \tilde{h}]$. Working by induction it can be shown that for all $\bar{d} \in L^\infty_{loc}(\Re^+; \Re^+)$ the solution of $\dot{x}(t) = a(x(t)) + \tilde{k}(x(\tau_i))$, $\tau_{i+1} = \tau_i + \tilde{h} \exp(-\bar{d}(\tau_i))$ starting at $x(0) = x_0 \in (0, 2)$ satisfies $x(t) \in \Theta$ and $|x(t)| \leq \exp(-t)|x_0|$, for all $t \geq 0$.

On the other hand, the solution of $\dot{x} = a(x) + \tilde{k}(x_0)$ starting at $x(0) = x_0 < 0$ satisfies $x_0 \leq x(t) \leq 0$ for all $t \geq 0$. Consequently, it holds that $\frac{d}{dt}|x(t)| = -a(-|x(t)|)$ for all $t \geq 0$ and Lemma 4.4 in [21] implies the existence of $\sigma \in KL$ with $\sigma(\sigma(s,t), \tau) = \sigma(s, t+\tau)$ for all $s, t, \tau \geq 0$ such that $|x(t)| \leq \sigma(|x_0|, t)$ for all $t \geq 0$. Working inductively, it can be shown that for all $\bar{d} \in L_{loc}^\infty(\Re^+; \Re^+)$ the solution of $\dot{x}(t) = a(x(t)) + \tilde{k}(x(\tau_i))$, $\tau_{i+1} = \tau_i + \tilde{h} \exp(-\bar{d}(\tau_i))$ starting at $x(0) = x_0 < 0$ satisfies $x(t) \in \Theta$ and $|x(t)| \leq \sigma(|x_0|, t)$, for all $t \geq 0$.

Therefore, hypothesis (P1) holds for system (4.7).

We next show that hypothesis (P3) holds as well. Consider the sets $\Omega_1 = \Theta$, $\Omega_j = (j-1, j]$ for $j = 2, 3, \ldots$. We will show that for all $j = 2, 3, \ldots$ and $r > 0$, the set $\bigcup_{i=1}^{j-1} \Omega_i \subseteq \Re$ is $r$-robustly reachable from the set $\Omega_j \subseteq \Re$ for system (4.7) with constant control. Notice that $\bigcup_{i=1}^{j-1} \Omega_i = \Theta$ for $j = 2$ and $\bigcup_{i=1}^{j-1} \Omega_i = (-\infty, j-1]$ for $j \geq 3$. Let $v_j = -1 - \max_{j-1 \leq s \leq j} a(s)$ and consider the solution of $\dot{x} = a(x) + v_j$ with initial condition $x(0) = x_0 \in \Omega_j$. As long as the solution exists, the following inequalities hold: $\dot{x} \geq v_j$ and $x(t) \leq x_0$. Consequently, it holds that $x_0 + t v_j \leq x(t) \leq x_0$. A simple contradiction argument shows that the solution exists for all $t \geq 0$ and satisfies $|x(t)| \leq \exp(t)(|x_0| + g(2|x_0|))$, where $g(s) = s + \max_{0 \leq x \leq s} a(x)$ (a function of class $K_\infty$). For $j \geq 3$, the fact that the set $\bigcup_{i=1}^{j-1} \Omega_i = (-\infty, j-1]$ is $r$-robustly reachable from the set $\Omega_j \subseteq \Re$ can be shown by following the procedure in the proof of Lemma 2.7 (with $V(x) = x$). For $j = 2$, the fact that the set $\bigcup_{i=1}^{j-1} \Omega_i = (-\infty, 2)$ is $r$-robustly reachable from the set $\Omega_j \subseteq \Re$ can be shown by the fact that the solution satisfies $x(t) < x_0$ for all $t > 0$.

Thus system (4.7) is robustly globally stabilizable by means of sampled-data control with positive sampling rate. A possible selection of the feedback is:

$$k(x) = \tilde{k}(x), \text{ for } x \in (-\infty, 2)$$
$$k(2) = -1 - \max_{1 \leq s \leq 2} a(s)$$
$$k(x) = -1 - \max_{j-1 \leq s \leq j} a(s), \text{ for } x \in (j-1, j], \ j \geq 3$$

where $\tilde{k}(x)$ is defined by (4.8) and $h = \frac{1}{L+1}$ where $L > 0$ is the Lipschitz constant that satisfies (4.9). ◁

It should be emphasized that the methodology used in Example 4.1 and Example 4.2, which is based on the main results of the present work (Theorem 3.1 and Theorem 3.2):

∗ provides a *simple formula* for a stabilizing sampled-data feedback,
∗ guarantees *global* asymptotic stability for the closed-loop system,
∗ guarantees *robustness* to perturbations of the sampling schedule,
∗ provides means to *determine the maximum allowable sampling period*,
∗ is *not limited* to special cases where the solution map is available,
∗ is *not limited* to special cases where the nonlinear term is homogeneous or globally Lipschitz

No other sampled-data feedback design methodology available in the literature can provide all the above features simultaneously.

# Appendix

**Proof of Lemma 3.4:** Lemma 4.4 in [21], guarantees the existence of a continuous function $\sigma$ of class $KL$, with $\sigma(s,0) = s$ for all $s \geq 0$ which satisfies $\frac{\partial}{\partial t}\sigma(s,t) = -\rho(\sigma(s,t))$ for all $s,t \geq 0$ with the following property: if $y:[t_0,t_1] \to \Re^+$ is an absolutely continuous function and $I \subset [t_0,t_1]$ a set of Lebesgue measure zero such that $\dot{y}(t)$ is defined on $[t_0,t_1]\setminus I$ and such that the following differential inequality holds for all $t \in [t_0,t_1]\setminus I$:

$$\dot{y}(t) \leq -\rho(y(t)) \tag{A1}$$

then the following estimate holds for all $t \in [t_0,t_1]$:

$$y(t) \leq \sigma(y(t_0), t-t_0) \tag{A2}$$

Actually, the statement of Lemma 4.4 in [21] does not guarantee that $\sigma$ is continuous or that $\sigma(s,0) = s$ for all $s \geq 0$, but a close look at the proof of Lemma 4.4 in [21] shows that this is the case when $\rho: \Re^+ \to \Re^+$ is a positive definite continuous function.

Let $x_0 \in \Theta$, $x_0 \neq 0 \in \Re^n$, $(d,\tilde{d}) \in L^\infty_{loc}(\Re^+;D) \times L^\infty_{loc}(\Re^+;\Re^+)$. The solution $x(t)$ of (2.4) with $k \equiv \tilde{k}$ exists locally and by virtue of (3.2) satisfies $x(t) \in \Theta$ and $M|x(t) - x_0| \leq |x(t)|$ for $t > 0$ sufficiently small. Let $\tau_1 = h\exp(-\tilde{d}(0))$ and notice that as long as the conditions $x(t) \in \Theta$ and $M|x(t) - x_0| \leq |x(t)|$ hold for $t \in [0,\tau_1)$ we have from (3.3) and (3.4):

$$\frac{d}{dt}V(x(t)) \leq -\rho(V(x(t))) \text{ and } \frac{d}{dt}\left(\frac{1}{2}|x(t)-x_0|^2\right) \leq L|x(t)-x_0|^2 + \gamma|x_0|^2, \text{ a.e.} \tag{A3}$$

Consequently, by virtue of (A1), (A2) and (A3) we obtain the following inequalities which hold as long as the conditions $x(t) \in \Theta$ and $M|x(t) - x_0| \leq |x(t)|$ hold for $t \in [0,\tau_1]$:

$$V(x(t)) \leq \sigma(V(x_0),t) \text{ and } |x(t)-x_0| \leq |x_0|\sqrt{\gamma\frac{\exp(2Lt)-1}{L}} \text{ for the case } L > 0 \tag{A4}$$

or

$$V(x(t)) \leq \sigma(V(x_0),t) \text{ and } |x(t)-x_0| \leq |x_0|\sqrt{2\gamma t} \text{ for the case } L = 0 \tag{A5}$$

By using (A4), (A5) in conjunction with the trivial inequality $|x_0| \leq |x(t)-x_0| + |x(t)|$, we obtain the following inequalities which hold as long as the conditions $x(t) \in \Theta$ and $M|x(t)-x_0| \leq |x(t)|$ hold for $t \in [0,\tau_1]$:

$$V(x(t)) \leq \sigma(V(x_0), t) \text{ and } |x(t) - x_0| \leq \frac{\sqrt{\gamma(\exp(2Lt)-1)}}{\sqrt{L} - \sqrt{\gamma(\exp(2Lt)-1)}} |x(t)| \text{ for the case } L > 0 \quad (A6)$$

or

$$V(x(t)) \leq \sigma(V(x_0), t) \text{ and } |x(t) - x_0| \leq \frac{\sqrt{2\gamma t}}{1 - \sqrt{2\gamma t}} |x(t)| \text{ for the case } L = 0 \quad (A7)$$

Notice that since $t \leq \tau_1 \leq h < \frac{1}{2L} \ln\left(1 + \frac{L}{\gamma}\left(\frac{1}{1+M}\right)^2\right)$, for the case $L > 0$ or $t \leq \tau_1 \leq h < \frac{1}{2\gamma}\left(\frac{1}{1+M}\right)^2$ for the case $L = 0$, we conclude that $M \frac{\sqrt{\gamma(\exp(2Lt)-1)}}{\sqrt{L} - \sqrt{\gamma(\exp(2Lt)-1)}} < 1$, $\sqrt{\gamma \frac{\exp(2Lt)-1}{L}} < 1$ for the case $L > 0$ or $M \frac{\sqrt{2\gamma t}}{1 - \sqrt{2\gamma t}} < 1$, $\sqrt{2\gamma t} < 1$ for the case $L = 0$. Consequently, combining the previous inequalities with (A4), (A5), (A6) and (A7), we obtain the following estimate for the solution $x(t)$ of (2.4) with $k \equiv \tilde{k}$, which holds as long as the conditions $x(t) \in \Theta$ and $M|x(t) - x_0| \leq |x(t)|$ hold for $t \in [0, \tau_1]$:

$$V(x(t)) \leq \sigma(V(x_0), t) \text{ and } M|x(t) - x_0| < |x(t)|$$

A simple contradiction argument shows that the above estimate holds for all $t \in [0, \tau_1]$. Notice that estimate $V(x(t)) \leq \sigma(V(x_0), t)$ holds for the case $x_0 = 0$ as well.

Using induction and the semigroup property for $\sigma$, we obtain for all non-negative integers $i$:

$$V(x(t)) \leq \sigma(V(x_0), t), \text{ for all } t \in [\tau_i, \tau_{i+1}] \quad (A8)$$

Since $V: \Re^n \to \Re^+$ is a positive definite, continuously differentiable and radially unbounded function, there exist functions $a_1, a_2 \in K_\infty$ such that $a_1(|x|) \leq V(x) \leq a_2(|x|)$, for all $x \in \Re^n$ (Lemma 3.5, page 138 in [16]). The conclusion of Lemma 3.2 is an immediate consequence of (A8), the previous inequality and the properties of KL functions. The proof is complete. ◁

**Proof of Fact in the proof of Theorem 3.1:** Let $(d, \tilde{d}) \in L^\infty_{loc}(\Re^+; D) \times L^\infty_{loc}(\Re^+; \Re^+)$ and $x(\tau_i) \in C_k$ with $k > 1$ (arbitrary). By virtue of definition (3.5) if follows that $x(\tau_i) \in \Omega_k$. Let $v_k \in U$ be the constant control that guarantees property (Q) of Definition 2.4 with $\Omega$ replaced by $\Omega_k$ and $A = \bigcup_{i=1}^{k-1} \Omega_i = B_{k-1}$. The solution of (2.4) on $[\tau_i, \tau_{i+1}]$ coincides with the solution of (1.1) with $u(t) \equiv v_k$, same initial conditions and corresponding to the same input $d \in L^\infty_{loc}(\Re^+; D)$.

Let $T(d, x(\tau_i)) \in [0, c + b(|x(\tau_i)|)]$ the time involved in property (Q) of Definition 2.4. Since $\lim_{p \to \infty} \tau_{i+p} = +\infty$ (see [12]), there exists integer $p \geq 0$ such that $\tau_{i+p+1} - \tau_i > T(d, x(\tau_i)) \geq \tau_{i+p} - \tau_i$. We claim that there exists non-negative integer $q \leq p+1$ such that $x(\tau_{i+q}) \in C_s$ for some $s < k$. Notice that an immediate consequence of the claim is that $\tau_{i+q} \leq \tau_{i+p+1} \leq \tau_{i+p} + r \leq \tau_i + T(d, x(\tau_i)) + r \leq \tau_i + c + b(|x(\tau_i)|) + r$.

By virtue of property (Q) of Definition 2.4 and since $\tau_{i+1} \leq \tau_i + r$, the solution of (2.4) exists for all $t \in [\tau_i, \tau_{i+1}]$ and satisfies $|x(t)| \leq a(|x(\tau_i)|)$, for all $t \in [\tau_i, \tau_{i+1}]$. In order to prove the above claim we distinguish the following cases:

a) $T(d, x(\tau_i)) > \tau_{i+1} - \tau_i$. In this case we have $x(\tau_{i+1}) \in \Omega_k \subseteq B_k$. Since $B_k = \bigcup_{s=1,\ldots,k} C_s$, it follows that there exists $s \in \{1,\ldots,k\}$ such that $x(\tau_{i+1}) \in C_s$.

b) $T(d, x(\tau_i)) \leq \tau_{i+1} - \tau_i$. Since $\tau_{i+1} - \tau_i \leq r$, in this case there exists $m \in \{1,...,k-1\}$ such that $x(\tau_{i+1}) \in \Omega_m \subseteq B_m$. Since $B_m = \bigcup_{s=1,...,m} C_s$, there exists $s < k$ such that $x(\tau_{i+1}) \in C_s$.

In every case we obtain the existence of $s \in \{1,...,k\}$ such that $x(\tau_{i+1}) \in C_s$. However, if $x(\tau_{i+1}) \in C_k \subseteq \Omega_k$, then $T(d, x(\tau_i)) > \tau_{i+1} - \tau_i$ and thus we can guarantee that property (Q) of Definition 2.4 holds with $T(d, x(\tau_{i+1})) = T(d, x(\tau_i)) - (\tau_{i+1} - \tau_i)$. Furthermore, since $\tau_{i+2} - \tau_{i+1} \leq r$, the solution of (2.4) exists for all $t \in [\tau_i, \tau_{i+2}]$ and satisfies $|x(t)| \leq a(|x(\tau_i)|)$, for all $t \in [\tau_i, \tau_{i+2}]$. By distinguishing cases (similarly as above), we conclude that there exists $s \in \{1,...,k\}$ such that $x(\tau_{i+2}) \in C_s$.

However, if $x(\tau_{i+2}) \in C_k$, then $T(d, x(\tau_{i+1})) > \tau_{i+2} - \tau_{i+1}$ and thus we can guarantee that property (Q) of Definition 2.4 holds with $T(d, x(\tau_{i+2})) = T(d, x(\tau_{i+1})) - (\tau_{i+2} - \tau_{i+1}) = T(d, x(\tau_i)) - (\tau_{i+2} - \tau_i)$. Furthermore, since $\tau_{i+3} - \tau_{i+2} \leq r$, the solution of (2.4) exists for all $t \in [\tau_i, \tau_{i+3}]$ and satisfies $|x(t)| \leq a(|x(\tau_i)|)$, for all $t \in [\tau_i, \tau_{i+3}]$.

Continuing in the same way, we conclude that there exists non-negative integer $q \leq p$ such that $x(\tau_{i+q}) \in C_s$ for some $m \in \{1,...,k-1\}$, because otherwise we would have $T(d, x(\tau_{i+p+1})) = T(d, x(\tau_i)) - (\tau_{i+p+1} - \tau_i) < 0$ (a contradiction). Moreover, the solution of (2.4) satisfies $|x(t)| \leq a(|x(\tau_i)|)$, for all $t \in [\tau_i, \tau_{i+q}]$.

The proof is complete. ◁